\documentclass{article}
\usepackage{mathrsfs}
\usepackage{amsmath}
\usepackage{amssymb, amsmath}
\usepackage{graphicx}
\catcode`\@=11 \@addtoreset{equation}{section}

\catcode`\@=12
\usepackage{color}

\newtheorem{Theorem}{Theorem}[section]
\newtheorem{Proposition}{Proposition}[section]
\newtheorem{Lemma}{Lemma}[section]
\newtheorem{Corollary}{Corollary}[section]
\newtheorem{Definition}{Definition}[section]

\newtheorem{Remark}{Remark}[section]

\newcommand{\newcom}{\newcommand}
\newcommand{\bTheorem}[1]{
\begin{Theorem} \label{T#1} }
\newcommand{\eT}{\end{Theorem}}

\newcommand{\bProposition}[1]{
\begin{Proposition} \label{P#1}}
\newcommand{\eP}{\end{Proposition}}

\newcommand{\bLemma}[1]{
\begin{Lemma} \label{L#1} }
\newcommand{\eL}{\end{Lemma}}

\newcommand{\bCorollary}[1]{
\begin{Corollary} \label{C#1} }
\newcommand{\eC}{\end{Corollary}}

\newcommand{\beq}{\begin{equation}}
\newcommand{\eeq}{\end{equation}}
\newcom{\ben}{\begin{eqnarray}}
\newcom{\een}{\end{eqnarray}}
\newcom{\beno}{\begin{eqnarray*}}
\newcom{\eeno}{\end{eqnarray*}}
\newcom{\bali}{\begin{aligned}}
\newcom{\eali}{\end{aligned}}

\newcommand{\bFormula}[1]{
\begin{equation} \label{#1}}
\newcommand{\eF}{\end{equation}}

\newcommand{\f}{\frac}
\newcommand{\Om}{\Omega}

\newcommand{\p}{\partial}

\newcommand{\vr}{\varrho}

\newcommand{\vu}{\vc{u}}

\newcommand{\vv}{\vc{v}}

\newcommand{\vc}[1]{{\boldsymbol #1}}
\newcommand{\Div}{{\rm div}}
\newcommand{\Grad}{\nabla}

\newcommand{\dx}{{\rm d} x}
\newcommand{\dy}{{\rm d} y}
\newcommand{\dt}{{\rm d} t }
\newcommand{\ds}{{\rm d} s}

\newcommand{\dxdt}{\dx\dt}

\font\F=msbm10 scaled 1000
\newcommand{\R}{\mbox{\F R}}

\newcommand\Cbox[2]{%
    \newbox\contentbox%
    \newbox\bkgdbox%
    \setbox\contentbox\hbox to \hsize{%
        \vtop{
            \kern\columnsep
            \hbox to \hsize{%
                \kern\columnsep%
                \advance\hsize by -2\columnsep%
                \setlength{\textwidth}{\hsize}%
                \vbox{
                    \parskip=\baselineskip
                    \parindent=0bp
                    #2
                }%
                \kern\columnsep%
            }%
            \kern\columnsep%
        }%
    }%
    \setbox\bkgdbox\vbox{
        \color{#1}
        \hrule width  \wd\contentbox %
               height \ht\contentbox %
               depth  \dp\contentbox
        \color{black}
    }%
    \wd\bkgdbox=0bp%
    \vbox{\hbox to \hsize{\box\bkgdbox\box\contentbox}}%
    \vskip\baselineskip%
}

\begin{document}


\title{\bf Global weak solutions and long time behavior for 1D compressible MHD equations without resistivity}

\author{Yang Li \\ Department of Mathematics, \\ Nanjing University, Nanjing 210093, China \\ lynjum@163.com \\
Yongzhong Sun \\ Department of Mathematics, \\ Nanjing University, Nanjing 210093, China \\ sunyz@nju.edu.cn}

\maketitle
{\bf Abstract: }{We study the initial-boundary value problem for 1D compressible MHD equations of viscous non-resistive fluids in the Lagrangian mass coordinates. Based on the estimates of upper and lower bounds of the density, weak solutions are constructed by approximation of global regular solutions, the existence of which has recently been obtained by Jiang and Zhang in \cite{JZW}. Uniqueness of weak solutions is also proved as a consequence of Lipschitz continuous dependence on the initial data. Furthermore, long time behavior for global solutions is investigated. Specifically, based on the uniform-in-time bounds of the density from above and below away from zero, together with the structure of the equations, we show the exponential decay rate in $L^2$- and $H^1$-norm respectively, with initial data of arbitrarily large.}


\section{Introduction}
The motion of conducting fluids is described by the system of  magnetohydrodynamics (MHD).  In Eulerian coordinates, a typical model for 3D compressible MHD fluids assumes the following form (see \cite{HC}):
\beq\label{a1}
\vr_t + \Div(\vr\vu)=0,
\eeq
\beq\label{a2}
(\vr\vu)_t+\Div(\vr\vu \otimes \vu)+\Grad{p}=\nu\Delta\vu + (\nu +\eta)\Grad{\Div \vu}
+(\Grad \times \vc{b})\times \vc{b},
\eeq
\beq\label{a3}
\vc{b}_t=\Grad \times (\vu \times \vc{b})-\lambda \Grad \times (\Grad \times \vc{b}),
\eeq
\beq\label{a4}
\Div \vc{b}=0.
\eeq
Here the unknown functions $\vr,\,\vu \in \R^3,\,p$ and $\vc{b}\in\R^3$ denote the density of fluid, the velocity, the pressure and the magnetic field, respectively. The viscosity coefficients $\nu$ and $\eta$ satisfy
\[
\nu>0 ,\,\,3\eta+2\nu \geq 0.
\]
Moreover, $\lambda\ge 0$ is the resistivity coefficient which represents the magnetic diffusion of the field $\vc{b}$. The compressible fluid is assumed to be isentropic, which means the pressure $p$ is prescribed through the following constitutive relation:
\beq\label{a5}
p(\vr)=A\vr^{\gamma},
\eeq
where $A$ is a positive constant and the adiabatic exponent $\gamma >1$.

Assuming that the resistivity coefficient $\lambda$ is positive, based on the pioneering work of P. L. Lions \cite{LS}, E. Feireisl et al. \cite{FNP}, Hu and Wang \cite{HW} obtained the global existence and large time behavior of global weak solutions to 3D equations of compressible MHD flows. B. Ducomet and E. Feireisl \cite{DF} proved the global weak solutions to the Navier-Stokes-Fourier system, coupled with the Maxwell equations with finite energy initial data.

In fact, the resistivity coefficient $\lambda$ is extremely small in practical models, and the fluid is often referred to as perfect conductor if $\lambda=0$. Thus, it is reasonable to consider the compressible isentropic MHD equations without resistivity, where (\ref{a3}) reads as
\[
\vc{b}_t=\Grad \times (\vu \times \vc{b}).
\]

Compared with the case of positive resistivity, mathematical investigations to (\ref{a1})-(\ref{a4}) with $\lambda=0$ are relatively few. Obviously, zero resistivity  introduces extra difficulty to build global solutions. The only known results on the multi-dimensional case is the recent work of Wu and Wu \cite{WW}, where the authors have established global well-posedness for the initial value problem of 2D compressible non-resistive MHD system with initial data close to the stationary solution $\vu_s=(0,0),\vc{b}_s=(1,0)$. This is an extension of early results of Lin et al. \cite{LXZ,XZ} for incompressible MHD to the compressible case. As to the incompressible MHD without resistivity, see also \cite{RWZ,ZT}.

In this article, we focus on the MHD equations without resistivity and restrict ourselves to the simplest one-dimensional case. By assuming
\[
\vr=\vr(x,t),\,\vu=(u(x,t),0,0),\,\vc{b}=(0,0,b(x,t)),
\]
where $x\in \R$ is the spatial variable, (\ref{a1})-(\ref{a4}) (with $\lambda=0$) are reduced to (see \cite{G})
\beq\label{a10}
\vr_t+(\vr u)_x=0,
\eeq
\beq\label{a11}
(\vr u)_t+(\vr u^2 +p(\vr)+\f{1}{2}b^2)_x=\mu u_{xx},
\eeq
\beq\label{a12}
b_t+(b u)_x=0,
\eeq
where the pressure $p$ satisfies (\ref{a5}) and $\mu=2\nu +\eta>0$.

Recently, Jiang and Zhang in \cite{JZW}  obtained the global well-posedness of strong solutions to the initial-boundary value problem for (\ref{a10})-(\ref{a12}) with initial data of arbitrary size, by making a full use of the effective viscous flux, the material derivative and the structure of the equations. See also Yu \cite{YU} for a similar result concerning the appearance of vacuum, but with more restriction on the initial magnetic field. We refer to \cite{FH,ZZ} for more results on 1D compressible heat-conductive MHD equations with vanishing resistivity.

It should be noted that if the resistivity coefficient $\lambda$ is included above, (\ref{a12}) becomes
\beq\label{a13}
b_t+(b u)_x=\lambda b_{xx}.
\eeq
There are many investigations for system (\ref{a10})-(\ref{a11}), (\ref{a13}). Kazhikhov and Smagulov in \cite{KS} announced the global well-posedness of strong solutions to the one-dimensional compressible, heat-conductive, viscous fluids with resistivity. Fan, Jiang and G. Nakamura in \cite{FJGN} obtained the existence, uniqueness and Lipschitz continuous dependence on the initial data of global weak solutions to a similar system.

As it is well-known, the classical method to handle one-dimensional models in fluid mechanics is the use of Lagrangian mass coordinates. To this end, we assume the fluid occupies the interval $[0,1]$. Let
\beq\label{a14}
y=\int_0^x\vr(\xi,t)d\xi,\,\,s=t.
\eeq
Then (\ref{a10})-(\ref{a12}) are reformulated as
\beq\label{a15}
\tau_t=u_x,
\eeq
\beq\label{a16}
u_t=\sigma_x,
\eeq
\beq\label{a17}
(b\tau)_t=0,
\eeq
with
\[
\tau :=\vr^{-1}
\]
the specific volume of the flow and the effective viscous flux
\beq\label{a18}
\sigma:=\mu\f{u_x}{\tau}-A\tau^{-\gamma}-\f{1}{2}b^2.
\eeq
Here for convenience, we still use $(x,t)$ instead of $(y,s)$ to denote the spatial and temporal variables. Without loss of generality, we assume the conserved total mass on $[0,1]$ is one unit. We then supplement system (\ref{a15})-(\ref{a18}) with the following initial and boundary conditions:
\beq\label{a19}
\tau(x,0)=\tau_0(x),u(x,0)=u_0(x),b(x,0)=b_0(x),x\in[0,1],
\eeq
\beq\label{a20}
u(0,t)=u(1,t)=0,\,\, t\in (0,\infty).
\eeq

The present paper is dedicated to the study of global weak solutions to the initial-boundary value problem (\ref{a15})-(\ref{a20}). Based on estimates of upper and lower bounds of the density, we first construct weak solutions by approximation of global regular solutions, the existence of which is guaranteed by Jiang and Zhang \cite{JZW}. Then we show the stability of weak solutions, that is, the Lipschitz continuous dependence on the initial data. The uniqueness of global weak solutions follows as a consequence of stability. In particular, similar to the results for one-dimensional Navier-Stokes(-Fourier) system, see \cite{ AZ2,CD, Hoff2, JZ, ZA1} among others, our results show that neither vacuum nor concentration can form in finite time for weak solutions. Furthermore,  based on the uniform-in-time bounds of the density from above and below away from zero, the exponential decay estimates of solutions are obtained in $L^2$- and $H^1$-norm respectively.

It should be noted that the stabilization for 1D compressible barotropic Navier-Stokes equations has been well-established since the work of Kanel \cite{KN1} and Kazhikhov \cite{KZ}. Extensions to more general barotropic case or the inclusion of external forces can be found in \cite{HBV1,HBV2,MY1,SV,SZ1,Z1}. The reader may consult \cite{AKS} for the stability of 1D Navier-Stokes-Fourier system in bounded domain. Also see \cite{JS1,KN2,LJ} for the case of unbounded domains.

Before giving the main results of this paper, we introduce the notations and functional spaces used throughout this paper. Denote $\Omega:=(0,1)$, $\Omega_t:=\Omega\times(0,t)$. Let $p\in[1,\infty]$, $k$ be a positive integer. We denote the usual Lebesgue space $L^p(\Omega)$  by $L^p$, with its norm $||\cdot||_{L^p}$; $H^k$ denotes the usual Sobolev space $H^k(\Omega)$, with its norm $||\cdot||_{H^k}$; $L^p(0,T;X)$ is the space of all strongly measurable, $p$th-power integrable functions from $(0,T)$ to $X$, with $X$ being some Banach space and its corresponding norm $||\cdot||_{L^p(0,T;X)}$. The Sobolev space $W^{1,p}(0,T;X)$ consists of all functions $\vv \in L^p(0,T;X)$ such that $\vv_t$ exists in the weak sense and belongs to $L^p(0,T;X)$. The Banach space $C([0,T];X)$ stands for all continuous functions from $[0,T]$ to $X$.

Concerning with the initial-boundary value problem (\ref{a15})-(\ref{a20}) for an isentropic, viscous and compressible flow, the first result of this paper is the existence of global weak solutions.
\begin{Theorem}\label{ly1}
Assume
\beq\label{a25}
\inf\limits_{x\in (0,1)}\tau_0(x) > 0, \,\tau_0\in L^{\infty}, \,u_0\in L^2,\, b_0\in L^{\infty}.
\eeq
Then there exists a weak solution $(\tau,u,b)$ to (\ref{a15})-(\ref{a20}) in the time interval $[0,T]$ for any fixed $T\in (0,\infty)$. Moreover, there exists a constant $C>0$, such that
\beq\label{a26}
C^{-1}\leq \tau(x,t) \leq C,\, |b(x,t)|\leq C,\text{ for a.e. }(x,t)\in \Omega_T,
\eeq
\beq\label{a127}
\|u\|_{L^{\infty}(0,T;L^2)}+\|u_x\|_{L^2(0,T;L^2)}+\|\tau_t\|_{L^2(0,T;L^2)}+\|b_t\|_{L^2(0,T;L^2)}\leq C.
\eeq
\end{Theorem}

Here and in the next theorem, the letter $C$ denotes a generic positive constant depending only on the parameters $A,\gamma,\mu$, the fixed time $T$ and the initial data. The definition of weak solutions will be given in the next section.

The next theorem concerns the stability of weak solutions obtained in Theorem \ref{ly1}.

\begin{Theorem}\label{ly2}
Let $(\tau,u,b)$ and $(\widetilde{\tau},\widetilde{u},\widetilde{b})$ be two weak solutions on $[0,1]\times[0,T]$ corresponding to the initial data $(\tau_0,u_0,b_0)$ and $(\widetilde{\tau}_0,\widetilde{u}_0,\widetilde{b}_0)$.
Then there exists a constant $C>0$ such that
\[
\|\tau-\widetilde{\tau}\|_{L^{\infty}(\Omega_T)}+\|u-\widetilde{u}\|_{L^{\infty}(0,T;L^2)}+
\|b-\widetilde{b}\|_{L^{\infty}(\Omega_T)}
+\|(u-\widetilde{u})_x\|_{L^2(0,T;L^2)}
\]
\beq\label{a30}
\leq C(\|\tau_0-\widetilde{\tau}_0\|_{L^{\infty}}+\|b_0-\widetilde{b}_0\|_{L^{\infty}}+
\|u_0-\widetilde{u}_0\|_{L^2}).
\eeq
\end{Theorem}

Obviously, Theorem \ref{ly2} in particular implies the uniqueness of weak solutions.
\begin{Corollary}\label{c1}
Under the assumptions of Theorem \ref{ly1}, there exists a unique global weak solution to the initial-boundary value problem (\ref{a15})-(\ref{a20}).
\end{Corollary}

The subsequent two theorems are associated with the long time behavior for global solutions to (\ref{a15})-(\ref{a20}).
\begin{Theorem}\label{sl1}
 Let the assumption (\ref{a25}) be satisfied and $(\tau,u,b)$ be the unique weak solution to (\ref{a15})-(\ref{a20}). Then there exist two positive constants $C_1,C_2$ which are independent of time, such that
\beq\label{a206}
\|(\tau-\tau_s)(t)\|_{L^2}+\|u(t)\|_{L^2}+\|(b-b_s)(t)\|_{L^2}\leq C_1\exp (-C_2 t),\text{ for any }t\geq 0.
\eeq
\end{Theorem}

Here, $(\tau_s,0,b_s)$ are the stationary solution to (\ref{a15})-(\ref{a20}) which will be introduced in Section \ref{s6}. Here and in the next theorem, we denote $C,C_i$ by generic positive constants depending only on the parameters of the system, the initial data and the stationary solution.

Given more regular initial data, we are able to strengthen the exponential decay rate of solutions in $H^1$-norm. To be more precise, we have
\begin{Theorem}\label{sl2}
Assume
\beq\label{a27}
\inf\limits_{x\in (0,1)}\tau_0(x) > 0, \,\tau_0\in W^{1,\infty}, \,u_0\in H^1_0,\, b_0\in W^{1,\infty}.
\eeq
Let $(\tau,u,b)$ be the unique strong solution to (\ref{a15})-(\ref{a20}). Then there exist two positive constants $C_3,C_4$, such that
\beq\label{a28}
\|(\tau-\tau_s)(t)\|_{H^1}+\|u(t)\|_{H^1}+\|(b-b_s)(t)\|_{H^1}\leq C_3\exp (-C_4 t),\text{ for any }t\geq 0.
\eeq
\end{Theorem}

The key point to obtain these results, especially Theorem \ref{ly1} and \ref{sl1} on the existence of global weak solution and its long time behavior, is the observation that under the Lagrangian formulation, the magnetic field ${b}$ is solved out as $b={b_0\tau_0}\tau^{-1}$. This observation results in the momentum equation a non-standard pressure law $p=p(x,\tau)$. The dependence of $p$ on the spatial variable $x$ makes it difficult to apply the traditional approaches for 1D isentropic Navier-Stokes equations such as in \cite{KZ}, especially for uniform pointwise estimates for the density. To overcome this difficulty, we have to modify the methods developed in \cite{AKS} to handle the full Navier-Stokes-Fourier system as well as in \cite{Z1} to treat a wider class of pressure laws. Moreover, it is also new for the large time behavior of the specific volume as well as the magnetic field in the sense that they approach to the nontrivial stationary solution $(\tau_s,b_s)$ determined by (\ref{wy20}) and (\ref{wy21}).

The rest of this paper is organized as follows. In Section \ref{s2} we recall the existence of global strong solution due to Jiang and Zhang \cite{JZW} under the framework of Lagrangian coordinates. In Section \ref{s3} we prove Theorem \ref{ly1} by approximation of strong solutions. The proof of Theorem \ref{ly2} is completed in Section \ref{s4} by modifying the ideas used in \cite{JZ}. The proof of Theorem \ref{sl1} and Theorem \ref{sl2} are finished in Section \ref{s6} by means of establishing the necessary uniform-in-time estimates.

\section{Preliminary Results}\label{s2}
To establish the existence of weak solution, we use approximation of strong solutions, the existence of which has been obtained in \cite{JZW} in the framework of Eulerian coordinates. It should be pointed out that for the initial-boundary value problem (\ref{a15})-(\ref{a20}), the global existence (and uniqueness) of strong solutions still holds in our case of Lagrangian formulation. Here, for completeness and later use, we just state this result and give a sketch of the proof. Throughout the present and the next two sections, the letter $C$ denotes a generic positive constant which is described after the statement of Theorem \ref{ly1} in the introduction.
\begin{Proposition}\label{js}
Assume that the initial data $(\tau_0,u_0,b_0)$ given in (\ref{a19}) satisfy
\beq\label{b1}
 \min_{0\leq x\leq 1}\tau_0(x)>0,\,(\tau_0,b_0)\in H^1,\, u_0\in H^1_0.
\eeq
Then there exists a unique strong solution $(\tau,u,b)$ in the time interval $[0,\infty)$ to the initial-boundary value problem (\ref{a15})-(\ref{a20}) such that
\beq\label{b3}
(\tau,b)\in L^\infty_{loc}(0,\infty;H^1),   (\tau_t,b_t)\in L^2_{loc}(0,\infty; L^2),
\eeq
\beq\label{b4}
u\in L^\infty_{loc}(0,\infty;H^1_0)\cap L^2_{loc}(0,\infty;H^2), \, u_t\in L^2_{loc}(0,\infty;L^2).
\eeq
Furthermore, for any fixed $0<T<\infty$, there exists a positive constant $C$ such that
\beq\label{b2}
C^{-1}\leq \tau(x,t) \leq C,\text{ for any }(x,t)\in[0,1]\times[0,T],
\eeq
\beq\label{b5}
\|\tau,b,u\|_{L^\infty(0,T;H^1)} + \|\tau_t,b_t,u_{xx},u_t\|_{L^2(\Omega_T)} \leq C.
\eeq
\end{Proposition}

The proof of this proposition is essentially based on global a priori estimates. We first give the standard energy estimates without proof.
\begin{Lemma}\label{wk1}
Let $(\tau,u,b)$ be a smooth solution to the initial-boundary value problem (\ref{a15})-(\ref{a20}) on $[0,1]\times[0,T]$. Then
\beq\label{wy1}
\int^1_0\tau(x,t)\dx=\int^1_0\tau_0(x)\dx=1,\text{ for any } t\in [0,T],
\eeq
and
\beq\label{wy2}
\sup_{0\leq t \leq T}\int^1_0\left(\f{1}{2}u^2+\f{A}{\gamma-1}\tau^{1-\gamma}+
\f{1}{2}b^2_0\tau^2_0\tau^{-1}\right)\dx+\mu\int_0^T\int_0^1\f{u_x^2}{\tau}\dx\ds \leq C.
\eeq
\end{Lemma}

The next lemma gives the upper and lower bounds of the specific volume, which is essential for the proof of Proposition \ref{js}. Here we modify the argument of Antontsev et. al., see \cite{AKS}.
\begin{Lemma}\label{wk2}
Let $(\tau,u,b)$ be a smooth solution to the initial-boundary value problem (\ref{a15})-(\ref{a20}) on $[0,1]\times[0,T]$. Then
\beq\label{wy3}
C^{-1}\leq \tau(x,t) \leq C,\text{ for any }(x,t)\in[0,1]\times[0,T].
\eeq
\end{Lemma}
{\bf Proof.} Note that (\ref{a16}) can be rewritten, using (\ref{a15}), (\ref{a17}) and (\ref{a18}), as
\[
u_t  =  \left(\mu\f{u_x}{\tau}-A\tau^{-\gamma}-\f{1}{2}b_0^2\tau_0^2\tau^{-2}\right)_x
\]
\[
     =  \left(\mu(\log\tau)_t-A\tau^{-\gamma}-\f{1}{2}b_0^2\tau_0^2\tau^{-2}\right)_x.
\]
By (\ref{wy1}) and the mean value theorem, for each $t\in [0,T]$, there exists $a(t)\in [0,1]$, such that $\tau(a(t),t)=1$. Integrating the above equation first over $(0,t)$ with respect to $t$, then over $(a(t),x)$ ($x$ is an arbitrarily fixed point in $[0,1]$) with respect to $x$, and then taking exponential on both sides of the resulting equation, we obtain the following representation of the specific volume $\tau$ as follows.
\beq\label{wy4}
Y(t)\tau(x,t)=\tau_0(x)B(x,t) \exp \left\{ \f{1}{\mu} \int_0^t \left( A\tau^{-\gamma}+\f{1}{2}b_0^2\tau_0^2\tau^{-2} \right) (x,s) \ds \right\},
\eeq
where
\[
B(x,t) := \exp\left(\f{1}{\mu}\int^x_{a(t)}u(\xi,t)-u_0(\xi) \, d\xi\right);
\]
\[
Y(t):=\tau_0(a(t))\exp\left\{\f{1}{\mu}\int_0^t\left( A \, \tau^{-\gamma}+\f{1}{2}b_0^2\tau_0^2\tau^{-2}\right)(a(t),s)\ds\right\}.
\]
By Lemma \ref{wk1} and Cauchy-Schwarz's inequality,
\beq\label{wy5}
C^{-1}\leq B(x,t) \leq C, \, C^{-1} \leq Y(t),\text{ for any }(x,t)\in[0,1]\times[0,T].
\eeq

We compute
\[
\frac{\partial}{\p t}\exp\left\{\f{1}{\mu}\int_0^t\left( A\tau^{-\gamma}+\f{1}{2}b_0^2\tau_0^2\tau^{-2}\right)(x,s)\ds\right\}
\]
\[ =  \f{1}{\mu}\left(A\tau^{-\gamma}+\f{1}{2}b_0^2\tau_0^2\tau^{-2}\right)\exp\left\{\f{1}{\mu}\int_0^t\left( A\tau^{-\gamma}+\f{1}{2}b_0^2\tau_0^2\tau^{-2}\right)(x,s)\ds\right\}
\]
\[ = \f{1}{\mu}\left(A\tau^{-\gamma}+\f{1}{2}b_0^2\tau_0^2\tau^{-2}\right)Y(t)\tau(x,t)\tau_0^{-1}(x)B^{-1}(x,t).
\]
Integrating the above equation over $(0,t)$ with respect to $t$ gives
\[
\exp\left\{\f{1}{\mu}\int_0^t\left( A\tau^{-\gamma}+\f{1}{2}b_0^2\tau_0^2\tau^{-2}\right)(x,s)\ds\right\}
\]
\[
=1+\f{1}{\mu}\int_0^t\left( A\tau^{1-\gamma}+\f{1}{2}b_0^2\tau_0^{2}\tau^{-1}\right)(x,s)Y(s)\tau_0^{-1}(x)B^{-1}(x,s)\ds.
\]
By substituting the above identity into (\ref{wy4}), we find
\[
Y(t)\tau(x,t)=\tau_0(x)B(x,t)
\]
\beq\label{wy6}
 \times \left\{1+\f{1}{\mu}\int_0^t\left( A\tau^{1-\gamma}+\f{1}{2}b_0^2\tau_0^{2}\tau^{-1}\right)(x,s)Y(s)\tau_0^{-1}(x)B^{-1}(x,s)\ds\right\}.
\eeq
Integrating (\ref{wy6}) over $(0,1)$ with respect to $x$ and by virtue of (\ref{wy1}), (\ref{wy2}) and (\ref{wy5}),
\[
Y(t)\leq C\left(1+\int_0^t Y(s)\ds\right),
\]
which together with Gronwall's inequality yields
\beq\label{wy7}
Y(t)\leq C,\text{ for any }t\in[0,T].
\eeq
Therefore, (\ref{wy4}), (\ref{wy5}) and (\ref{wy7}) imply
\beq\label{wy8}
C^{-1}\leq \tau(x,t), \text{ for any }(x,t)\in[0,1]\times[0,T].
\eeq
The upper bound of the specific volume $\tau$ follows immediately from (\ref{wy5})-(\ref{wy8}). This completes the proof of Lemma \ref{wk2}.
\begin{Remark}\label{re1}
We note that in \cite{JZW}, to deal with the vanishing resistivity problem in Eulerian coordinates, the authors have to use a different approach to show the boundedness of the density from above and below by making a full use of the effective viscous flux, the material derivative and the structure of the equations. In particular, in their proof the lower boundedness of the density follows from the boundedness of the magnetic field, while in our case the boundedness of the magnetic field follows directly from that of the specific volume obtained in Lemma \ref{wk2}.
\end{Remark}

Once we have Lemma \ref{wk1} and Lemma \ref{wk2} at hand, it remains to derive the higher order energy estimates for the specific volume $\tau$, the magnetic field $b$, and the velocity field $u$. We list the higher order energy estimates with detailed proof omitted here,  see \cite{AKS,JZW}.
\begin{Lemma}\label{wk3}
Let $(\tau,u,b)$ be a smooth solution to the initial-boundary value problem (\ref{a15})-(\ref{a20}) on $[0,1]\times[0,T]$. Then there exists a constant $C>0$ such that
\beq
\|\tau_t,\tau_x, b_t, b_x, u_x\|_{L^\infty(0,T;L^2)} \leq C,
\eeq
\beq
\|u_{xx},u_t,\tau_{tx},b_{tx}\|_{L^2(\Omega_T)} \leq C.
\eeq
\end{Lemma}
Based on these a priori estimates, the global existence of strong solutions to the initial-boundary value problem (\ref{a15})-(\ref{a20}) can be proved in a standard way.

Finally we introduce the definition of weak solution to the MHD system (\ref{a15})-(\ref{a17}).
\begin{Definition}\label{df1}
We say that  $(\tau,u,b)$ is a weak solution to MHD system (\ref{a15})-(\ref{a17}) on $[0,1]\times [0,T]$ with boundary condition (\ref{a20}) and initial data $(\tau_0,u_0,b_0)$ satisfying
\[
\inf\limits_{x\in (0,1)}\tau_0(x) > 0, \, \tau_0, b_0\in L^{\infty}, \, u_0\in L^2,
\]
provided that
\[
\tau \in W^{1,2}(0,T;L^2), \,u\in L^\infty(0,T;L^2)\cap L^2(0,T;H^1_0), \,b\in W^{1,2}(0,T;L^2),
\]
\[
\inf\limits_{(x,t)\in \Omega_T}\tau(x,t)> 0,\, \tau, b \in L^{\infty}(\Omega_T),
\]
\[
\tau_t=u_x, \,  b=b_0 \tau_0 \tau^{-1}  \text{ a.e. in }\Omega_T,
\]
\[
\tau(x,0)=\tau_0(x) \text{ for a.e. } x \in \Om,
\]
and that for any test function $\psi \in C^{\infty}_c(\Om \times [0,T))$, the following integral identity holds:
\[
\int_0^T \int_{\Omega} u\psi_t-\left(\mu\f{u_x}{\tau}-A\tau^{-\gamma}-\f{1}{2}b^2\right)\psi_x \,\dxdt+\int_{\Om}u_0\psi(x,0)\dx=0.
\]
\end{Definition}

\section{Existence of weak solutions}\label{s3}
In this section, to prove Theorem \ref{ly1}, we first obtain a strong solution to the initial-boundary value problem (\ref{a15})-(\ref{a20}) by regularizing the initial data and then show the existence of global weak solutions.

Under the assumptions of initial data in (\ref{a25}), we construct a sequence $(\tau_0^{\epsilon},u_0^{\epsilon},b_0^{\epsilon})$, by regularizing the initial data, such that
\[
b_0^{\epsilon},\,\tau_0^{\epsilon}\in C^2([0,1]),\, u_0^{\epsilon}\in C^2_c((0,1)),
\]
\[
C^{-1}\leq \tau_0^{\epsilon}\leq C,\, (\tau_0^{\epsilon},u_0^{\epsilon},b_0^{\epsilon}) \,
\rightarrow(\tau_0,u_0,b_0)\text{ strongly in } L^2\text{ as }\epsilon\rightarrow 0^{+},
\]
\[
\|b_0^{\epsilon}\|_{L^{\infty}}\leq \|b_0\|_{L^{\infty}}.
\]

Now we consider the initial-boundary value problem (\ref{a15})-(\ref{a20}) with $(\tau_0,u_0,b_0)$ replaced by the approximate initial data $(\tau_0^{\epsilon},u_0^{\epsilon},b_0^{\epsilon})$. It follows from Proposition \ref{js} that there exists a unique global strong solution $(\tau^{\epsilon},u^{\epsilon},b^{\epsilon})$ such that
\[
0<\tau^{\epsilon}<\infty \text{ in } [0,1]\times[0,T],
\]
\[
u^{\epsilon}\in L^\infty(0,T;L^2),u^{\epsilon}_x\in L^2(0,T;L^2),
\]
\[
\tau^{\epsilon}_t\in L^2(0,T;L^2),\,b^{\epsilon}_t\in L^2(0,T;L^2).
\]

It should be pointed out that a careful review of Lemmas \ref{wk1}-\ref{wk2} shows that the approximate solutions $(\tau^{\epsilon},u^{\epsilon},b^{\epsilon})$  have the following uniform-in-$\epsilon$ bounds:
\beq\label{wj1}
C^{-1}\leq \tau^{\epsilon}(x,t) \leq C, \text{ for any }(x,t)\in [0,1]\times[0,T],
\eeq
\beq\label{wj2}
\|u^{\epsilon}\|_{L^\infty(0,T;L^2)}+\|u^{\epsilon}_x\|_{L^2(0,T;L^2)}\leq C,
\eeq
\beq\label{wj3}
\|\tau^{\epsilon}_t\|_{L^2(0,T;L^2)}+\|b^{\epsilon}_t\|_{L^2(0,T;L^2)}\leq C.
\eeq

In order to pass to the limits to obtain the existence of weak solutions to (\ref{a15})-(\ref{a20}), we have to show that the specific volume $\tau$ exists as a strong limit of $\tau^{\epsilon}$, due to the nonlinearity of the system. For this purpose, we give the following crucial lemma. Let $\Delta_{h}w(x):=w(x+h)-w(x)$, which is the difference of $w$ with respect to $x$.
\begin{Lemma}\label{yd1}
For any $0<h<1$, there holds
\[
\|\Delta_h\tau^{\epsilon}\|_{L^\infty(0,T;L^2)}\leq C(\|\Delta_h \tau_0\|_{L^2}+\|\Delta_h b_0\|_{L^2}+h).
\]
\end{Lemma}
{\bf Proof.} Set
\[
 \sigma^{\epsilon}:=\mu \f{u^{\epsilon}_x}{\tau^{\epsilon}}-A(\tau^{\epsilon})^{-\gamma}-\f{1}{2}(b^{\epsilon})^2,\,
 a^{\epsilon}_0:=b_0^{\epsilon}\tau_0^{\epsilon}.
 \]
Thus $(\tau^{\epsilon},u^{\epsilon},b^{\epsilon})$ satisfies the following system:
\beq\label{wj4}
\tau^{\epsilon}_t=u^{\epsilon}_x,
\eeq
\beq\label{wj5}
u^{\epsilon}_t= \sigma^{\epsilon}_x,
\eeq
\beq\label{wj6}
(b^{\epsilon}\tau^{\epsilon})_t=0.
\eeq
Note that (\ref{wj4})-(\ref{wj6}) together give us
\beq\label{wj7}
\tau^{\epsilon}_t=\f{\sigma^{\epsilon}}{\mu}\tau^{\epsilon}+\f{A}{\mu}(\tau^{\epsilon})^{1-\gamma}+
\f{(a^{\epsilon}_0)^2}{2\mu\tau^{\epsilon}}.
\eeq
Multiplying (\ref{wj7}) by $\exp\left(-\f{1}{\mu}\int_0^t\sigma^{\epsilon}(x,s)\ds\right) $ and integrating the resulting equation over $(0,t)$ with respect to $t$ yields
\[
\exp\left(-\f{1}{\mu}\int_0^t\sigma^{\epsilon}(x,s)\ds\right) \tau^{\epsilon}
\]
\[
=\tau^{\epsilon}_0+\f{1}{\mu}
\int_0^t\exp\left(-\f{1}{\mu}\int_0^{\xi}\sigma^{\epsilon}(x,s)\ds\right)
\left(A(\tau^{\epsilon})^{1-\gamma}+ \f{(a^{\epsilon}_0)^2}{2\tau^{\epsilon}}\right)(x,\xi)d\xi.
\]
Hence we have
\[
\tau^{\epsilon} = \left[\tau^{\epsilon}_0 +\f{1}{\mu}
                   \int_0^t\exp\left(-\f{1}{\mu}\int_0^{\xi}\sigma^{\epsilon}(x,s)\ds\right)
                   \left(A(\tau^{\epsilon})^{1-\gamma}+ \f{(a^{\epsilon}_0)^2}{2\tau^{\epsilon}} \right)(x,\xi)d\xi\right]
\]
\beq\label{wj8}
               \times  \exp\left(\f{1}{\mu}\int_0^t\sigma^{\epsilon}(x,s)\ds\right).
\eeq

By defining
\[
B^{\epsilon}(x,t):=\exp\left(\f{1}{\mu}\int_0^t\sigma^{\epsilon}(x,s)\ds\right),
\]
and recalling (\ref{wj1}), one sees
\beq\label{wj9}
C^{-1}\leq B^{\epsilon}(x,t) \leq C, \text{ for any }(x,t)\in[0,1]\times[0,T].
\eeq
Consequently, (\ref{wj8}) reads
\[
\tau^{\epsilon}=B^{\epsilon}\left[\tau^{\epsilon}_0+\f{1}{\mu}
\int_0^t (B^{\epsilon})^{-1}(x,\xi)
\left(A(\tau^{\epsilon})^{1-\gamma}+ \f{(a^{\epsilon}_0)^2}{2\tau^{\epsilon}} \right)(x,\xi)d\xi\right],
\]
and direct computation shows that
\[
\Delta_h\tau^{\epsilon}(x,t)= \Delta_h B^{\epsilon}(x,t)\left[\tau^{\epsilon}_0(x+h)+\f{1}{\mu}
                             \int_0^t (B^{\epsilon})^{-1}
                             \left(A(\tau^{\epsilon})^{1-\gamma}
                             + \f{(a^{\epsilon}_0)^2}{2\tau^{\epsilon}}\right)(x+h,\xi)d\xi\right]
\]
\[
                             +  B^{\epsilon}(x,t)\left[\Delta_h\tau^{\epsilon}_0+\f{1}{\mu}\int_0^t (B^{\epsilon})^{-1}(x+h,\xi)\Delta_h\left(A(\tau^{\epsilon})^{1-\gamma}
                             +  \f{(a^{\epsilon}_0)^2}{2\tau^{\epsilon}}\right)(x,\xi) d\xi\right]
\]
\[
                             -  B^{\epsilon}(x,t)\left[ \f{1}{\mu} \int_0^t (B^{\epsilon})^{-1}(x+h,\xi)(B^{\epsilon})^{-1}(x,\xi)
                             \left(A(\tau^{\epsilon})^{1-\gamma}
                             +  \f{(a^{\epsilon}_0)^2}{2\tau^{\epsilon}}\right)(x,\xi)\Delta_h B^{\epsilon}(x,\xi) d\xi\right],
\]
which, by (\ref{wj1}) and (\ref{wj9}), implies
\beq
\|\Delta_h\tau^{\epsilon}(x,t)\|_{L^2} \leq  C(\|\Delta_hB^{\epsilon}(x,t)\|_{L^2}+
                                       \|\Delta_h\tau^{\epsilon}_0\|_{L^2})
\eeq
\beq
                                       + C\int_0^t\|\Delta_hB^{\epsilon}(x,\xi)\|_{L^2}
                                       +\|\Delta_h\tau^{\epsilon}(x,\xi)\|_{L^2}
                                       +\|\Delta_h a^{\epsilon}_0\|_{L^2} d\xi
\eeq
\[
 \leq  C\left(h\|u^{\epsilon}-u_0^{\epsilon}\|_{L^\infty(0,T;L^2)}+
 \|\Delta_h\tau^{\epsilon}_0\|_{L^2}
  + \|\Delta_h b^{\epsilon}_0\|_{L^2}
  +\int_0^t\|\Delta_h\tau^{\epsilon}(x,\xi)\|_{L^2} d\xi\right)
\]
\beq\label{wj10}
                                       \leq C \left(h+\|\Delta_h\tau_0\|_{L^2}+\|\Delta_h b_0\|_{L^2}+ \int_0^t\|\Delta_h\tau^{\epsilon}(x,\xi)\|_{L^2} d\xi\right).
\eeq
An application of Gronwall's inequality to (\ref{wj10}) yields
\[
\|\Delta_h\tau^{\epsilon}\|_{L^\infty(0,T;L^2)}\leq C(\|\Delta_h \tau_0\|_{L^2}+\|\Delta_h b_0\|_{L^2}+h),
\]
thus completing the proof of Lemma \ref{yd1}.

Note that (\ref{wj1})-(\ref{wj3}) allow us to extract a subsequence of $(\tau^{\epsilon},u^{\epsilon},b^{\epsilon})$, still denoted by $(\tau^{\epsilon},u^{\epsilon},b^{\epsilon})$, such that as $\epsilon\rightarrow0^{+}$, the following weakly or weakly-$\star$ convergences hold:
\beq\label{wj11}
\tau^{\epsilon}\rightarrow\tau \text{ weakly}-\star \text{ in }L^{\infty}(0,T;L^{\infty}),
\eeq
\beq\label{wj12}
u^{\epsilon}\rightarrow u  \text{ weakly}-\star \text{ in }L^{\infty}(0,T;L^2),
\eeq
\beq\label{wj13}
(\tau^{\epsilon}_t,u^{\epsilon}_x)\rightarrow(\tau_t,u_x) \text{ weakly in }L^2(0,T;L^2).
\eeq
In addition, for the limit functions $(\tau,u)$,  we have
\beq\label{wj14}
C^{-1}\leq \tau(x,t) \leq C , \text{ for a.e. }(x,t)\in\Om_T,
\eeq
\beq\label{wj15}
\|u\|_{L^{\infty}(0,T;L^2)}+\|u_x\|_{L^2(0,T;L^2)}+\|\tau_t\|_{L^2(0,T;L^2)}\leq C.
\eeq
By (\ref{wj3}) and Lemma \ref{yd1}, we deduce that for any $0<h<1$, $0<s<T$, there holds
\beq\label{wj16}
\|\tau^{\epsilon}(\cdot+h,\cdot+s)-\tau^{\epsilon}\|_{L^{\infty}(0,T-s;L^2)}
\leq C(\|\Delta_h \tau_0\|_{L^2}+\|\Delta_h b_0\|_{L^2}+h+s^{\f{1}{2}}).
\eeq
Recalling the criterion of compactness of sets in $L^2(0,T;L^2)$ and invoking (\ref{wj11}), (\ref{wj16}) implies
\beq\label{wj17}
\tau^{\epsilon}\rightarrow \tau \text{ strongly in } L^2(0,T;L^2) \text{ as }\epsilon\rightarrow0^{+}.
\eeq
By means of defining
\[
b:=b_0\tau_0\tau^{-1},
\]
one checks easily, by virtue of (\ref{wj1}), (\ref{wj13}), (\ref{wj14}) and (\ref{wj17}), that
\beq\label{wj18}
b^{\epsilon}\rightarrow b \text{ strongly in } L^2(0,T;L^2)\text{ as }\epsilon\rightarrow0^{+},
\eeq
\beq\label{wj19}
b^{\epsilon}_t\rightarrow b_t\text{ weakly in } L^2(0,T;L^2)\text{ as }\epsilon\rightarrow0^{+},
\eeq
\beq\label{wj20}
\|b_t\|_{L^2(0,T;L^2)}\leq C.
\eeq

Based on Lemma \ref{yd1} and the analysis of weak convergence given above, we are now ready to give the proof of Theorem \ref{ly1}.

We multiply (\ref{wj5}) by any $\phi \in C^{\infty}_c((0,1) \times [0,T))$, then integrate over $\Om_T$, and perform an integration by parts. Letting $\epsilon\rightarrow0^{+}$, taking (\ref{wj11})-(\ref{wj13}), (\ref{wj17})-(\ref{wj19}) into account, we find that $(\tau,u,b)$ obtained is a global weak solution to the initial-boundary value problem (\ref{a15})-(\ref{a20}), by gathering the results for $(\tau,u,b)$ derived above. Moreover, the estimates (\ref{a26}) and (\ref{a127}) follow from (\ref{wj14}), (\ref{wj15}) and (\ref{wj20}). The proof of Theorem \ref{ly1} is therefore complete.

\section{Uniqueness of weak solutions}\label{s4}
In this section, we prove Theorem \ref{ly2} by modifying the arguments used in \cite{AZ2,JZ}. The proof is based on the following three lemmas.
\begin{Lemma}\label{lf1}
Let the assumptions of Theorem \ref{ly2} be satisfied. Then the following representations are valid in $\Om_T$:
\[
\tau(x,t) = \exp\left(\f{1}{\mu}\int_0^t\sigma(x,s)\ds\right)
\]
\beq\label{kb1}
          \times  \left[ \tau_0+\int_0^t\exp\left(-\f{1}{\mu}\int_0^{\xi}\sigma(x,s)\ds\right)
          \left(\f{A}{\mu}\tau^{1-\gamma}+\f{1}{2\mu}b_0^2\tau_0^2\tau^{-1}\right)(x,\xi) d\xi \right],
\eeq
and
\beq\label{kb2}
\int_0^t\sigma(x,s)\ds=(J_{\Omega}(u-u_0))(x,t)+\int_0^t<\sigma(\cdot,s)> \ds,
\eeq
where the linear operator $J_{\Omega}$ is defined by
\[
J_{\Omega}w(x):=\int_0^x w(\xi)d\xi-<\int_0^x w(\xi)d\xi>,\,\, <w>:=\int_0^1w(x)\dx.
\]
\end{Lemma}
{\bf Proof.} Obviously, (\ref{a15})-(\ref{a18}), and Theorem \ref{ly1} imply the following relations:
\beq\label{kb3}
\tau_t=\f{\sigma}{\mu}\tau+\f{A}{\mu}\tau^{1-\gamma}+\f{1}{2\mu}b_0^2\tau_0^2\tau^{-1},
\eeq
\beq\label{kb4}
\left(\int_0^t\sigma(x,s)\ds\right)_x=u-u_0.
\eeq
Multiplying (\ref{kb3}) by $\exp\left(-\f{1}{\mu}\int_0^t\sigma(x,s)\ds\right)$ and integrating the resulting equation over $(0,t)$ with respect to $t$ gives (\ref{kb1}). In addition, applying the operator $J_{\Omega}$ to (\ref{kb4}) yields (\ref{kb2}) immediately.

Before stating the next lemma, for simplicity, we introduce the notations below.
\[
(\Delta \tau,\Delta u, \Delta b):=(\tau-\widetilde{\tau},u-\widetilde{u},b-\widetilde{b}),
\]
\[
(\Delta \tau_0,\Delta u_0, \Delta b_0):=(\tau_0-\widetilde{\tau_0},u_0-\widetilde{u_0},b_0-\widetilde{b_0}),
\]
\[
\Delta\sigma:=\sigma-\widetilde{\sigma},\,\,
\widetilde{\sigma}:=\mu\f{\widetilde{u}_x}{\widetilde{\tau}}-A(\widetilde{\tau})^{-\gamma}
-\f{1}{2}(\widetilde{b_0})^2(\widetilde{\tau_0})^2(\widetilde{\tau})^{-2},
\]
\[
g:=\exp\left(\f{1}{\mu}\int_0^t\sigma(x,s)\ds\right),\,\,
\widetilde{g}:=\exp\left(\f{1}{\mu}\int_0^t\widetilde{\sigma}(x,s)\ds\right),
\]
\[
K:=\f{A}{\mu}\tau^{1-\gamma}+\f{1}{2\mu}b_0^2\tau_0^2\tau^{-1},\,\,
\widetilde{K}:=\f{A}{\mu}(\widetilde{\tau})^{1-\gamma}+\f{1}{2\mu}(\widetilde{b_0})^2(\widetilde{\tau_0})^2
(\widetilde{\tau})^{-1},
\]
\[
\widetilde{\vr}:=(\widetilde{\tau})^{-1},\,\,\Delta\vr:=\vr-\widetilde{\vr}.
\]

Then our essential lemma with respect to the supremum norm of $\Delta\tau$ reads as follows.
\begin{Lemma}\label{lf2}
Let the assumptions of Theorem \ref{ly2} be fulfilled. Then for any  $t\in (0, T]$,
\[
\|\Delta\tau\|_{L^{\infty}(\Omega_t)} \leq C(\|\Delta \tau_0\|_{L^{\infty}}+\|\Delta b_0\|
                                      _{L^{\infty}}+\|\Delta u_0\|_{L^2}
\]
\beq\label{kb5}
+\|\Delta u\|_{L^{\infty}(0,t;L^2)}
+\|(\Delta u)_x\|_{L^2(0,t;L^2}).
\eeq
\end{Lemma}
{\bf Proof.} It follows from (\ref{a26}) that
\beq\label{kb6}
C^{-1}\leq g,\,\,\widetilde{g} \leq C.
\eeq
Direct computation, by (\ref{kb1}), shows that
\beq\label{kb7}
\Delta\tau  = g\left[\Delta \tau_0+\int_0^t K\left(\f{1}{g}-\f{1}{\widetilde{g}}\right)+\f{K-\widetilde{K}}
           {\widetilde{g}}d\xi\right]
           + (g-\widetilde{g})\left(\widetilde{\tau_0}+\int_0^t\f{\widetilde{K}}{\widetilde{g}}d\xi\right).
\eeq
Using (\ref{kb6}) and (\ref{a26}), we estimate
\[
|\Delta\tau|\leq C\left(|\Delta \tau_0|+\int_0^t\left|\int_0^{\xi}\Delta\sigma\ds\right|+|\Delta\tau|+|\Delta b_0|+|\Delta \tau_0| d\xi\right)+C\left|\int_0^t\Delta\sigma\ds\right|,
\]
which means
\beq\label{kb8}
|\Delta\tau|\leq C\left(|\Delta \tau_0|+|\Delta b_0|+\int_0^t\left|\int_0^{\xi}\Delta\sigma\ds\right|+|\Delta\tau| d\xi\right)
+C\left|\int_0^t\Delta\sigma\ds\right|.
\eeq
Obviously, (\ref{kb8}) yields the bound
\beq\label{kb9}
\|\Delta\tau(\cdot,t)\|_{L^{\infty}}\leq C\left(\|\Delta \tau_0\|_{L^{\infty}}
                                    +\|\Delta b_0\|_{L^{\infty}}+\left\|\int_0^{\xi}\Delta\sigma \ds\right\|_{L^{\infty}
                                    (\Omega_t)}
                                    +\int_0^t\|\Delta\tau(\cdot,\xi)\|_{L^{\infty}} d\xi \right).
\eeq
By virtue of (\ref{kb2}), we estimate the third term on the right hand side of (\ref{kb9}) in the following manner.
\beq\label{kb10}
\left\|\int_0^{\xi}\Delta\sigma \ds\right\|_{L^{\infty}(\Omega_t)}  \leq  \|J_{\Omega}\Delta u_0\|_{L^{\infty}}
                                                        +\|J_{\Omega}\Delta u\|_
                                                        {L^{\infty}(\Omega_t)}
                                                        + \|\Delta\sigma\|_{L^1(\Omega_t)}.
\eeq
It is easy to see
\beq\label{kb11}
\Delta\sigma=\mu\tau^{-1}(\Delta u)_x+\varphi,
\eeq
where
\[
\varphi:=\mu(\Delta \vr)\widetilde{u}_x-\mu\f{K}{\tau}+\mu\f{\widetilde{K}}{\widetilde{\tau}}.
\]
Thus, by invoking (\ref{a26}), we arrive at
\[
|\Delta\sigma|\leq C |(\Delta u)_x|+|\varphi|,
\]
\beq\label{kb12}
|\varphi|\leq C |\Delta\tau|(|\widetilde{u}_x|+1)+C (|\Delta b_0|+|\Delta \tau_0|).
\eeq
Consequently,
\beq\label{kb13}
\|\Delta\sigma\|_{L^1(\Omega_t)}\leq C\left(\|(\Delta u)_x\|_{L^1(\Omega_t)}+\|\Delta b_0\|_{L^{\infty}}+\|\Delta \tau_0\|_{L^{\infty}}+\int_0^t \zeta(s)\|\Delta\tau(\cdot,s)\|_{L^{\infty}}\ds\right),
\eeq
where
\[
\zeta(t):=\|\widetilde{u}_x(\cdot,t)\|_{L^2}+1.
\]
In accordance with (\ref{a127}), valid is
\beq\label{kb14}
\|\zeta\|_{L^2(0,T)}\leq C.
\eeq
By (\ref{kb10}) and (\ref{kb13}), (\ref{kb9}) is further estimated as follows.
\[
\|\Delta\tau(\cdot,t)\|_{L^{\infty}}\leq C(\|\Delta \tau_0\|_{L^{\infty}}
                                    +\|\Delta b_0\|_{L^{\infty}}+\|J_{\Omega}\Delta u_0\|_{L^{\infty}}
                                    +\|J_{\Omega}\Delta u\|_ {L^{\infty}(\Omega_t)})
\]
\beq\label{kb15}
                                    +C\left(\|(\Delta u)_x\|_{L^1(\Omega_t)}+
                                    \int_0^t \overline{\zeta}(s)\|\Delta\tau(\cdot,s)\|_{L^{\infty}}\ds\right),
\eeq
where
\[
\overline{\zeta}(t):=\zeta(t)+1.
\]
Applying again Gronwall's inequality to (\ref{kb15}) yields
\[
\|\Delta\tau\|_{L^{\infty}(\Omega_t)}  \leq  C(\|\Delta \tau_0\|_{L^{\infty}}
                                      +\|\Delta b_0\|_{L^{\infty}}+\|J_{\Omega}\Delta u_0\|_
                                      {L^{\infty}}
\]
\[
                                       +  \|J_{\Omega}\Delta u\|_ {L^{\infty}(\Omega_t)}
                                      +\|(\Delta u)_x\|_{L^1(\Omega_t)}),
\]
which combined with (\ref{kb2}), Cauchy-Schwarz's inequality implies that
\[
\|\Delta\tau\|_{L^{\infty}(\Omega_t)} \leq  C(\|\Delta \tau_0\|_{L^{\infty}}+\|\Delta b_0\|
                                      _{L^{\infty}}+\|\Delta u_0\|_{L^2}
\]
\[
                                       +  \|\Delta u\|_{L^{\infty}(0,t;L^2)}
                                      +\|(\Delta u)_x\|_{L^2(0,t;L^2}).
\]
This completes the proof of Lemma \ref{lf2}.

The next lemma concerns the energy estimate of $\Delta u$.
\begin{Lemma}\label{lf3}
Let the hypotheses of Theorem \ref{ly2} be satisfied. Then for any $t\in (0,T]$,
\[
\|\Delta u\|_{L^{\infty}(0,t;L^2)}+\|(\Delta u)_x\|_{L^2(0,t;L^2)} \leq C(\|\Delta u_0\|_{L^2} +
                                                                 \|\Delta b_0\|_{L^{\infty}}
                                                                 + \|\Delta \tau_0\|_{L^{\infty}}
\]
\beq\label{kb16}
                                                                  + \|\zeta\|\Delta\tau
                                                                 (\cdot,s)\|_{L^{\infty}} \|_{L^2(0,t)}).
\eeq
Here $\zeta(t)=\|\widetilde{u}_x(\cdot,t)\|_{L^2}+1$.
\end{Lemma}

{\bf Proof.} By (\ref{a16}) and (\ref{kb11}), we have
\beq\label{kb17}
(\Delta u)_t=[\mu \tau^{-1}(\Delta u)_x+\varphi]_x.
\eeq
In terms of multiplying (\ref{kb17}) by $\Delta u$ and integrating the resulting equation over $\Omega_t$, we obtain after integration by parts that
\beq\label{kb18}
\|\Delta u\|_{L^{\infty}(0,t;L^2)}+\|(\Delta u)_x\|_{L^2(0,t;L^2)}\leq C(\|\Delta u_0\|_{L^2}+\|\varphi\|_{L^2(\Omega_t)}),
\eeq
where (\ref{a26}) and Cauchy-Schwarz's inequality have been invoked.
We conclude readily, by virtue of (\ref{kb12}), that
\beq\label{kb19}
\|\varphi\|_{L^2(\Omega_t)}\leq C(\|\zeta\|\Delta\tau(\cdot,s)\|_{L^{\infty}} \|_{L^2(0,t)}+
\|\Delta b_0\|_{L^{\infty}}+ \|\Delta \tau_0\|_{L^{\infty}}).
\eeq
Thus, Lemma \ref{lf3} is proved by substituting (\ref{kb19}) into (\ref{kb18}).

Based on the previous lemmas, we now give the proof of Theorem \ref{ly2}.

Multiplying (\ref{kb5}) by $\f{1}{2C}$ and adding the resulting inequality to (\ref{kb16}) gives rise to
\[
\|\Delta\tau(\cdot,t)\|_{L^{\infty}}+\|\Delta u\|_{L^{\infty}(0,t;L^2)}+\|(\Delta u)_x\|_{L^2(0,t;L^2)}
\]
\[
\leq C(\|\Delta b_0\|_{L^{\infty}}+ \|\Delta \tau_0\|_{L^{\infty}}+\|\Delta u_0\|_{L^2}
   +\|\zeta\|\Delta\tau(\cdot,s)\|_{L^{\infty}} \|_{L^2(0,t)}),
\]
which implies
\[
(\|\Delta\tau(\cdot,t)\|_{L^{\infty}}+\|\Delta u\|_{L^{\infty}(0,t;L^2)}+\|(\Delta u)_x\|_{L^2(0,t;L^2)})^2
\]
\beq\label{kb20}
\leq C\left((\|\Delta b_0\|_{L^{\infty}}+ \|\Delta \tau_0\|_{L^{\infty}}+\|\Delta u_0\|_{L^2})^2 + \int_0^t \zeta^2(s)\|\Delta\tau(\cdot,s)\|_{L^{\infty}}^2 \ds\right).
\eeq
An application of Gronwall's inequality to (\ref{kb20}) gives
\[
\int_0^t \zeta^2(s)\|\Delta\tau(\cdot,s)\|_{L^{\infty}}^2 \ds\leq C(\|\Delta b_0\|_{L^{\infty}}+ \|\Delta \tau_0\|_{L^{\infty}}+\|\Delta u_0\|_{L^2})^2.
\]
As a consequence, we obtain
\[
\|\Delta\tau\|_{L^{\infty}(\Omega_t)}+\|\Delta u\|_{L^{\infty}(0,t;L^2)}+\|(\Delta u)_x\|_{L^2(0,t;L^2)}
\]
\beq\label{kb21}
\leq C(\|\Delta b_0\|_{L^{\infty}}+ \|\Delta \tau_0\|_{L^{\infty}}+\|\Delta u_0\|_{L^2}).
\eeq
In addition, it follows from Definition \ref{df1} and (\ref{a26}) that
\beq\label{kb22}
\|\Delta b\|_{L^{\infty}(\Omega_t)}\leq C(\|\Delta \tau_0\|_{L^{\infty}}+\|\Delta b_0\|_{L^{\infty}}
+\|\Delta\tau\|_{L^{\infty}(\Omega_t)}).
\eeq
Thus (\ref{a30}) is verified if we multiply (\ref{kb22}) by $\f{1}{2C}$ and add the resulting inequality to
(\ref{kb21}).

\begin{Remark}\label{re2}
In fact, as the classical results on one-dimensional compressible Navier-Stokes-Fourier system, the weak solution $(\tau,u,b)$ obtained in Theorem \ref{ly1} satisfies
\[
\tau\in C([0,T];L^{\infty}),\,u\in C([0,T];L^2),\,b\in C([0,T];L^{\infty}),
\]
for any fixed $0<T<\infty$, see \cite{AZ2,ZA1}.
\end{Remark}

\section{Large time behavior}\label{s6}
The crucial step to the proof of Theorem \ref{sl1}, \ref{sl2}  lies in obtaining the uniform-in-time bounds of the density from above and below away from zero.  To this end, we first notice that the energy estimates given in Lemma \ref{wk1} are uniform with respect to time. For the sake of convenience, we rewrite it as follows.
\begin{Lemma}\label{le1}
Let $(\tau,u,b)$ be the unique weak solution to (\ref{a15})-(\ref{a20}) under the assumption (\ref{a25}). Then
\beq\label{m1}
\int^1_0\tau(x,t)\dx=\int^1_0\tau_0(x)\dx=1,\text{ for any } t\in [0,\infty),
\eeq
and
\beq\label{m2}
\sup_{0\leq t < \infty}\int^1_0\left(\f{1}{2}u^2+\f{A}{\gamma-1}\tau^{1-\gamma}+
\f{1}{2}b^2_0\tau^2_0\tau^{-1}\right)\dx+\mu\int_0^{\infty}\int_0^1\f{u_x^2}{\tau}\dx\ds \leq C.
\eeq
\end{Lemma}
During this section, the letter $C,C_i$ denote generic positive constants indepenent of the time.
Following Zlotnik \cite{Z1}, we first consider the boundary value problem with a parameter $t\geq 0$ as follows.
\beq\label{m3}
(\rho w_x)_x=f,\,\, x\in (0,1),\,\,w|_{x=0,1}=0.
\eeq
Here $\rho=\rho(t,x)>0$ and $f$ are given functions in $\Omega\times(0,\infty)$ and $w$ is the unknown function. Let $\eta=\rho^{-1}$ and $v$ satisfy $\eta_t=v_x$. Denote $\Lambda f:=\rho w_x$. We report the following results on $\Lambda$ from \cite{Z1}.
\begin{Lemma}\label{le2}
There holds
\beq\label{m4}
(\Lambda f^1)(x,t)=-\int_x^1 f^1(\xi,t) d\xi+\int_0^1\eta(y,t)\int_y^1f^1(\xi,t)d\xi\dy,
\eeq
\beq\label{m5}
(\Lambda f^2_x)(x,t)=f^2(x,t)-\int_0^1\eta(y,t)f^2(y,t)\dy,
\eeq
\beq\label{m6}
(\Lambda f^3)_t (x,t)=(\Lambda f^3_t)(x,t)+\int_0^1 v(y,t)f^3(y,t)\dy,
\eeq
\beq\label{m7}
\|\Lambda f^4\|_{L^{\infty}}\leq 2\|f^4\|_{L^1},
\eeq
where $f^1(\cdot,t),f^2(\cdot,t),f^4(\cdot,t),(\eta f^2)(\cdot,t)\in L^1$ for any $t\geq 0$ and $f^3,f^3_t,vf^3 \in L^1(\Om\times (0,T))$ for any $T \in (0,\infty)$.
\end{Lemma}

Based on Lemmas \ref{le1}-\ref{le2}, we can obtain the uniform-in-time bounds of the density from above and below away from zero, which plays a crucial role in deriving exponential decay estimates.
\begin{Lemma}\label{le3}
Let $(\tau,u,b)$ be the unique weak solution to (\ref{a15})-(\ref{a20}) under the assumption (\ref{a25}). Then
\beq\label{m8}
C^{-1}\leq \tau(x,t)\leq C, \text{ for any } (x,t)\in [0,1]\times [0,\infty).
\eeq
\end{Lemma}
{\bf Proof. } By setting
\[
P(x,\tau):=A\tau^{-\gamma}+\f{1}{2}b_0^2\tau_0^2\tau^{-2},
\]
we rewrite (\ref{a16}) as
\beq\label{wy9}
(\varrho u_x)_x=\f{1}{\mu}\left(u_t+P(x,\tau)_x\right).
\eeq
In view of Lemma \ref{le2} and (\ref{a20}),
\[
\varrho u_x=\f{1}{\mu}\left((\Lambda u)_t-\int_0^1u^2\dx+P(x,\tau)-\int_0^1\tau P(x,\tau)\dx\right).
\]
Thus, by (\ref{a15}), we see
\beq\label{wy10}
(\log \tau)_t=\f{1}{\mu}\left((\Lambda u)_t-\int_0^1u^2\dx+P(x,\tau)-\int_0^1\tau P(x,\tau)\dx\right).
\eeq

On the one hand, using (\ref{m2}), we get
\[
\int_0^1\tau P(x,\tau)\dx \leq C_1,
\]
and there exists $C_2>0$ such that
\beq\label{wy17}
P(x,\tau)> C_1,\text{ if }0<\tau<C_2.
\eeq
Now we fix $x\in [0,1]$ and set
\[
\underline{\tau_0}(x):=min\{\tau_0(x),C_2\}.
\]
If there exists $t_2\in (0,\infty)$ such that
\[
\tau(x,t_2)<\underline{\tau_0}(x),
\]
then, due to the continuity of $\tau(x,t)$ with respect to $t$ (see Remark \ref{re2}), there exists $t_1\in [0,t_2)$ such that
\beq\label{wy11}
\tau(x,t_1)=\underline{\tau_0}(x),\,\,\tau(x,t)<\underline{\tau_0}(x),\text{ for }t\in (t_1,t_2].
\eeq
Integrating (\ref{wy10}) both sides over $(t_1,t_2)$ with respect to $t$ yields
\[
\mu \log \tau(x,t_2)-\mu \log \tau(x,t_1)=(\Lambda u)(x,t_2)-(\Lambda u)(x,t_1)
\]
\beq\label{wy12}
-\int_{t_1}^{t_2}\int_0^1u^2 \dx\ds +\int_{t_1}^{t_2}\left(P(x,\tau)-\int_0^1 \tau P(x,\tau)\dx\right)\ds.
\eeq
Taking advantage of (\ref{m2}) and (\ref{m7}), one easily finds
\beq\label{wy13}
\|(\Lambda u )(t)\|_{L^{\infty}}\leq C, \text{ for any }t\in [0,\infty);
\eeq
while the third term on the right-hand side of (\ref{wy12}) can be estimated by
\[
\int_{t_1}^{t_2}\int_0^1u^2 \dx\ds \leq \int_{t_1}^{t_2}\|u\|_{L^{\infty}}^2\ds \leq \int_{t_1}^{t_2} \|u_x\|^2_{L^1}\ds
\]
\beq\label{wy14}
\leq \int_{t_1}^{t_2}\left(\int_0^1\f{u_x^2}{\tau}\dx\int_0^1 \tau\dx\right)\ds \leq C,
\eeq
where H\"{o}lder's inequality and Lemma \ref{le1} have been used. As a consequence, by gathering (\ref{wy17}), (\ref{wy11}), (\ref{wy13}) and (\ref{wy14}), we deduce from (\ref{wy12}) that
\beq\label{wy15}
\tau(x,t_2)\geq \underline{\tau_0}(x)\exp \left(-\f{C}{\mu}\right).
\eeq

On the other hand, by employing (\ref{m1}) and Jensen's inequality, there holds
\[
\int_0^1 \tau P(x,\tau)\dx \geq A \int_0^1 \tau ^{1-\gamma}\dx\geq A \left(\int_0^1\tau \dx\right)^{1-\gamma}\geq A,
\]
and there exists $C_3>0$ such that
\beq\label{wy18}
P(x,\tau)< A,\text{ if }\tau >C_3.
\eeq
Now we fix $x\in [0,1]$ and set
\[
\overline{\tau_0}(x):=max\{\tau_0(x),C_3\}.
\]
If there exists $t_2\in (0,\infty)$ such that
\[
\tau(x,t_2)> \overline{\tau_0}(x),
\]
then, due to the continuity of $\tau(x,t)$ with respect to $t$, there exists $t_1\in [0,t_2)$ such that
\beq\label{wy16}
\tau(x,t_1)=\overline{\tau_0}(x),\,\,\tau(x,t)>\overline{\tau_0}(x),\text{ for }t\in (t_1,t_2].
\eeq
Therefore, similar to the derivation of (\ref{wy15}), it follows from (\ref{wy13}), (\ref{wy14}), (\ref{wy18}) and (\ref{wy16}) that
\beq\label{wy19}
\tau(x,t_2)\leq  \overline{\tau_0}(x)\exp \left(\f{C}{\mu}\right).
\eeq
This completes the proof of Lemma \ref{le3} by combining (\ref{wy15}) with (\ref{wy19}).

Before turning to the proof of Theorem \ref{sl1}, we give the unique stationary solution of (\ref{a15})-(\ref{a20}) denoted by $(\tau_s,u_s,b_s)$, obeying
\beq\label{wy20}
\int_0^1\tau_s(x)\dx=\int_0^1\tau_0(x)=1,\,\,u_s=0,
\eeq
\beq\label{wy21}
A\tau_s^{-\gamma}+\f{1}{2}b_s^{2}=C_0,\,\,b_s=b_0\tau_0\tau_s^{-1},
\eeq
where the constant $C_0$ is determined by the normalized condition $\int_0^1\tau_s(x) \dx = 1$. It is obvious that $\tau_s$ is upper and lower bounded, i.e., there exists a constant $C>0$ such that
\beq\label{wy22}
C^{-1}\leq \tau_s \leq C.
\eeq
Furthermore, with the regularity class (\ref{a27}) imposed on the initial data, we have
\beq\label{wy23}
\|(\tau_s)_x\|_{L^{\infty}}\leq C.
\eeq

\subsection{Exponential decay in $L^2$-norm}

With Lemmas \ref{le1}-\ref{le3} at hand, we are now in a position to give the proof of Theorem \ref{sl1}. The proof is essentially based on the energy method by modifying the idea used in \cite{MY1,SZ1,Z1}, .

Firstly, owing to (\ref{wy21}), we rewrite (\ref{a16}) as
\beq\label{y1}
u_t+(P(x,\tau)-P(x,\tau_s))_x=(\mu \varrho u_x)_x.
\eeq
Multiplying (\ref{y1}) both sides by $u$ and integrating the resulting equation over $(0,1)$ with respect to $x$,
\beq\label{y0}
\f{d}{dt}\int_0^1 \left( \f{1}{2}u^2 +\f{A}{\gamma-1}\tau^{-\gamma+1}+A \tau_s^{-\gamma}\tau +\f{1}{2}b_0^2\tau_0^2\tau^{-1}+\f{1}{2}b_0^2\tau_0^2\tau_s^{-2}\tau\right) \dx
+\int_0^1 \mu \varrho u_x^2 \dx=0,
\eeq
where (\ref{a15}) is used.
Denote
\[
\Phi_1(\tau,\tau_s):=\f{A}{\gamma-1}\tau^{-\gamma+1}+A \tau_s^{-\gamma}\tau-\f{A\gamma}{\gamma-1}\tau_s^{-\gamma+1};
\]
\[
\Phi_2(x,\tau,\tau_s):=\f{1}{2}b_0^2\tau_0^2\left(\tau^{-1}+\tau_s^{-2}\tau-2\tau_s^{-1}\right).
\]
Then (\ref{y0}) is equivalent to
\beq\label{y2}
\f{d}{dt}\int_0^1  \left(\f{1}{2}u^2 +\Phi_1(\tau,\tau_s)+\Phi_2(x,\tau,\tau_s)\right)\dx+\int_0^1 \mu \varrho u_x^2 \dx=0,
\eeq
Note that $\Phi_1(\tau,\tau_s)$ can be written as
\[
\Phi_1(\tau,\tau_s)=A\tau_s^{-\gamma+1}G\left(\f{\tau}{\tau_s}\right),
\]
where
\[
G(z)=z-\f{\gamma}{\gamma-1}+\f{1}{\gamma-1}z^{-\gamma+1}.
\]
It follows that that
\[
G(1)=G'(1)=0,\,\,G^{\prime\prime}(z)>0, \text{ if }z >0.
\]
As a consequence, by invoking (\ref{m8}) and (\ref{wy22}), we conclude that
\beq\label{y3}
C^{-1}(\tau-\tau_s)^2\leq \Phi_1(\tau,\tau_s)\leq C (\tau-\tau_s)^2.
\eeq
Similarly it holds that
\beq\label{y4}
0 \leq \Phi_2(x,\tau,\tau_s) \leq C (\tau-\tau_s)^2.
\eeq

For a positive parameter $\varepsilon$, we multiply (\ref{y1}) both sides by $\varepsilon \int_0^x (\tau-\tau_s) d\xi$ and integrate the resulting equation over $(0,1)$ with respect to $x$ to find
\[
\f{d}{dt}\int_0^1 \varepsilon u \mathcal{J}(\tau-\tau_s)\dx-\varepsilon \int_0^1 \left(P(x,\tau)-P(x,\tau_s)\right)(\tau-\tau_s)\dx
\]
\beq\label{y5}
-\varepsilon \int_0^1 u^2 \dx +\varepsilon \int_0^1 \mu \varrho u_x (\tau-\tau_s)\dx=0,
\eeq
where, for simplicity, we have set
\[
\mathcal{J}(\tau-\tau_s):=\int_0^x (\tau-\tau_s) d\xi.
\]
Adding (\ref{y5}) to (\ref{y2}) yields
\[
\f{d}{dt}\int_0^1  \left(\f{1}{2}u^2 +\Phi_1(\tau,\tau_s)+\Phi_2(x,\tau,\tau_s)+\varepsilon u \mathcal{J}(\tau-\tau_s)\right)\dx
\]
\[
+\int_0^1 \mu \varrho u_x^2 \dx-\varepsilon \int_0^1 \left(P(x,\tau)-P(x,\tau_s)\right)(\tau-\tau_s)\dx
\]
\beq\label{y6}
=\varepsilon \int_0^1 u^2 \dx-\varepsilon \int_0^1 \mu \varrho u_x (\tau-\tau_s)\dx.
\eeq
Obviously, (\ref{m8}) and (\ref{wy22}) imply
\beq\label{y7}
-\varepsilon \int_0^1 \left(P(x,\tau)-P(x,\tau_s)\right)(\tau-\tau_s)\dx \geq C_1 \varepsilon \int_0^1 (\tau-\tau_s)^2 \dx.
\eeq
Similarly, by Cauchy-Schwarz's inequality and (\ref{m8}), we see
\beq\label{y8}
\left|\varepsilon \int_0^1 \mu \varrho u_x (\tau-\tau_s)\dx\right| \leq \f{C_2 \varepsilon}{2C_1}\int_0^1 \mu \varrho u_x^2 \dx +  \f{C_1 \varepsilon}{2} \int_0^1 (\tau-\tau_s)^2 \dx;
\eeq
\beq\label{y9}
C_3 \int_0^1 \mu \varrho u_x^2 \dx \geq  \int_0^1u^2\dx.
\eeq
In view of (\ref{y7})-(\ref{y9}), it follows from (\ref{y6}) that
\[
\f{d}{dt}\int_0^1  \left(\f{1}{2}u^2 +\Phi_1(\tau,\tau_s)+\Phi_2(x,\tau,\tau_s)+\varepsilon u \mathcal{J}(\tau-\tau_s)\right)\dx
\]
\beq\label{y10}
+\f{C_1 \varepsilon}{2} \int_0^1 (\tau-\tau_s)^2 \dx+\left[\left(1-\f{C_2 \varepsilon}{2C_1}\right)-C_3\varepsilon \right]\int_0^1 \mu \varrho u_x^2 \dx \leq 0.
\eeq
An application of Cauchy-Schwarz's inequality again shows
\beq\label{y11}
\left|\int_0^1\varepsilon u \mathcal{J}(\tau-\tau_s)\dx  \right| \leq \f{\varepsilon}{2}\int_0^1 u^2 \dx +\f{\varepsilon}{2}\int_0^1 (\tau-\tau_s)^2\dx.
\eeq
Therefore, defining
\[
\mathcal{E}:=\int_0^1  \left(\f{1}{2}u^2 +\Phi_1(\tau,\tau_s)+\Phi_2(x,\tau,\tau_s)+\varepsilon u \mathcal{J}(\tau-\tau_s)\right)\dx,
\]
and gathering (\ref{y3}), (\ref{y4}) and (\ref{y11}), after choosing $\varepsilon$ to be a sufficiently small constant, we arrive at
\beq\label{y12}
C_4^{-1} \left(\int_0^1 (\tau-\tau_s)^2 \dx+\int_0^1 u^2 \dx\right)\leq \mathcal{E}\leq C_4 \left(\int_0^1 (\tau-\tau_s)^2 \dx+\int_0^1 u^2 \dx\right).
\eeq

Finally, combining (\ref{y10}) with (\ref{y12}) gives
\[
\f{d}{dt}\mathcal{E}+C_5 \left(\int_0^1 (\tau-\tau_s)^2 \dx+\int_0^1 u^2 \dx\right)\leq 0,
\]
from which one obtains the decay estimate after using (\ref{y12}) and integration
\beq\label{y13}
\|(\tau-\tau_s)(t)\|_{L^2}+\|u(t)\|_{L^2}\leq C_6\exp (-C_7 t),\text{ for any }t\geq 0.
\eeq
Due to (\ref{m8}) and (\ref{wy22}), there holds
\beq\label{y14}
\|(b-b_s)(t)\|_{L^2}\leq C \|(\tau-\tau_s)(t)\|_{L^2}.
\eeq
This completes the proof of Theorem \ref{sl1} by adding (\ref{y14}) to (\ref{y13}).

At this stage, we intend to give an interesting remark concerning a special case of Theorem \ref{sl1}.
\begin{Remark}\label{re3}
Suppose, in addition to (\ref{a25}), if the absolute value of the ratio between the initial magnetic field and density is a positive constant, then the stationary magnetic field will be a piecewise constant. As a simple example, assume there exists a positive constant $\theta$ such that
\begin{equation}
{b_0(x)}{\tau_0(x)}=
\left\{
\begin{aligned}
\theta,\text{ if }x\in \left(0,\f{1}{2}\right), \\
-{\theta}, \text{ if }x\in \left(\f{1}{2},1\right). \\
\end{aligned}
\right.\nonumber
\end{equation}
Then, in accordance with (\ref{wy20})-(\ref{wy21}), the stationary solution exactly takes
\begin{equation*}
(\tau_s,u_s)=(1,0), \, b_s(x)=
\left\{
\begin{aligned}
\theta,\text{ if }x\in \left(0,\f{1}{2}\right), \\
-\theta,\text{ if }x\in \left(\f{1}{2}, 1\right). \\
\end{aligned}
\right.
\end{equation*}

\end{Remark}

\subsection{Exponential decay in $H^1$-norm}\label{s7}
Inspired by the method introduced in \cite{SV,SZ1,Z1}, we give the proof of Theorem \ref{sl2} in this section. To this end, we need the uniform-in-time bound of the density and the velocity in $H^1$-norm. Specifically, we have
\begin{Lemma}\label{le4}
Let $(\tau,u,b)$ be the unique strong solution to (\ref{a15})-(\ref{a20}), with the hypotheses of Theorem \ref{sl2} be satisfied. Then
\beq\label{z1}
\|\tau_x\|_{L^2}+\|u_x\|_{L^2} \leq C.
\eeq
\end{Lemma}
{\bf Proof. } Denoting
\[
F:=u-\mu (\log \tau)_x, \,
a(x):=b_0^2(x)\tau_0^2(x),
\]
we rewrite (\ref{a16}), by means of (\ref{a15}), as
\beq\label{z2}
F_t +(P(x,\tau))_x=0.
\eeq
By using the fact that
\[
(P(x,\tau))_x = \frac{\partial P}{\partial x} + \frac{\partial P}{\partial\tau} \tau_x
= \frac{\partial P}{\partial x} - \frac{\partial P}{\partial\tau} \frac{\tau}{\mu}F - \frac{\partial P}{\partial\tau} \frac{\tau}{\mu}u,
\]
we multiply (\ref{z2}) both sides by $F$ and integrate the resulting equation over $(0,1)$ to find
\[
\f{d}{dt} \int_0^1\left(\f{1}{2}F^2\right)\dx+\f{1}{\mu}\int_0^1 (A\gamma \tau^{-\gamma}+a(x)\tau^{-2})F^2\dx
\]
\[
=-\int_0^1 \left(\f12a'(x)\tau^{-2}F\right)\dx+\f{1}{\mu}\int_0^1 (A\gamma \tau^{-\gamma}+a(x)\tau^{-2})F u \dx,
\]
the right-hand side of which can be estimated by
\[
\left|-\int_0^1 \left(\f12a'(x)\tau^{-2}F\right)\dx\right| \leq C \|a'\|_{L^2}\|F\|_{L^2}\leq \delta_1 \|F\|_{L^2}^2+C_{\delta_1}\|a'\|_{L^2}^2;
\]
\[
\left|\f{1}{\mu}\int_0^1 (A\gamma \tau^{-\gamma}+a(x)\tau^{-2})F u \dx \right| \leq C \|u\|_{L^2}\|F\|_{L^2}\leq \delta_2 \|F\|_{L^2}^2+C_{\delta_2}\|u\|_{L^2}^2,
\]
where (\ref{m8}) and Cauchy-Schwarz's inequality have been used. Hence, by choosing $\delta_1,\delta_2$ to be sufficiently small and invoking (\ref{m2}), we obtain
\beq\label{z3}
\f{d}{dt} \int_0^1 F^2 \dx +C_1 \int_0^1 F^2 \dx \leq C_2.
\eeq
It follows that
\[
\int_0^1 F^2 \dx \leq C,
\]
which particularly implies, recalling (\ref{m2}) and (\ref{m8}), that
\beq\label{z4}
\int_0^1\tau_x^2\dx \leq C.
\eeq

To proceed, we write (\ref{a16}) as
\[
u_t+(P(x,\tau))_x=\mu \varrho_x u_x +\mu \varrho u_{xx},
\]
followed by multiplying both sides by $u_{xx}$, integrating over $(0,1)$ with respect to $x$. Then after integration by parts we see
\[
\f{1}{2}\f{d}{dt} \int_0^1u_x^2 \dx+\int_0^1 \mu \varrho u^2_{xx}\dx=\int_0^1 ( P(x,\tau))_x u_{xx}\dx-\int_0^1 \mu \varrho_x u_x u_{xx}\dx
\]
\beq\label{z5}
\leq \|u_{xx}\|_{L^2}\|( P(x,\tau))_x\|_{L^2}+C \|\tau_x\|_{L^2}\|u_{xx}\|_{L^2}\|u_{x}\|_{L^{\infty}}.
\eeq
Furthermore, owing to Gagliardo-Nirenberg inequality, one has
\beq\label{z6}
\|u_{x}\|_{L^{\infty}}\leq \delta_3 \|u_{xx}\|_{L^2}+C_{\delta_3}\|u\|_{L^2}.
\eeq
As a consequence, due to (\ref{m8}), (\ref{z4}), (\ref{z6}) and Cauchy-Schwarz's inequality, after choosing $\delta_3$ to be sufficiently small, we conclude from (\ref{z5}) that
\beq\label{z7}
\f{d}{dt} \int_0^1u_x^2 \dx+C_3\int_0^1 u_{xx}^2\dx \leq C_4 \left(\int_0^1u^2\dx+\int_0^1( P(x,\tau))_x^2\dx\right).
\eeq
Obviously, (\ref{m2}), (\ref{m8}) and (\ref{z4}) together lead to
\[
 C_4 \left(\int_0^1u^2\dx+\int_0^1( P(x,\tau))_x^2\dx\right) \leq C_5.
\]
In addition, since $\int_0^1 u_x \dx = 0$,
\[
\int_0^1 u_{x}^2 \dx \leq \int_0^1 u_{xx}^2 \dx.
\]
Thus, we strengthen (\ref{z7}) as
\[
\f{d}{dt} \int_0^1u_x^2 \dx+C_6\int_0^1u_x^2 \dx \leq C_7.
\]
Hence
\beq\label{z8}
\int_0^1 u_{x}^2 \dx \leq C.
\eeq
This completes the proof of Lemma \ref{le4} by adding (\ref{z8}) to (\ref{z4}).

Based on the previous lemmas, we are now in a position to prove Theorem \ref{sl2}.
Recall that (\ref{a16}) is equivalent to
\beq\label{z9}
u_t +(P(x,\tau)-P(x,\tau_s))_x=\mu \left(\log \left(\f{\tau}{\tau_s}\right)\right)_{xt}.
\eeq
We multiply (\ref{z9}) both sides by $\left(\log \left(\f{\tau}{\tau_s}\right)\right)_x$ and integrate the resulting equation over $(0,1)$ with respect to $x$, to infer that
\[
\f{d}{dt}\left\{\f{\mu}{2} \int_0^1 \left(\log \left(\f{\tau}{\tau_s}\right)\right)_x^2 \dx -\int_0^1 u\left(\log \left(\f{\tau}{\tau_s}\right)\right)_x \dx\right\}
\]
\beq\label{z10}
-\int_0^1 (P(x,\tau)-P(x,\tau_s))_x \left(\log \left(\f{\tau}{\tau_s}\right)\right)_x\dx=\int_0^1 \varrho u_x^2 \dx.
\eeq
To proceed,  we write the second term on the left-hand side of (\ref{z10}) as follows.
\[
-\int_0^1 (P(x,\tau)-P(x,\tau_s))_x \left(\log \left(\f{\tau}{\tau_s}\right)\right)_x\dx=\mathcal{R}_1+\mathcal{R}_2+\mathcal{R}_3,
\]
where
\[
\mathcal{R}_1:=A \gamma \int_0^1 [\tau^{-\gamma-1}\tau_x -\tau_s^{-\gamma-1}(\tau_s)_x ][\tau^{-1}\tau_x -\tau_s^{-1}(\tau_s)_x]\dx;
\]
\[
\mathcal{R}_2:=\int_0^1 [a\tau^{-3}\tau_x-a\tau_s^{-3}(\tau_s)_x][\tau^{-1}\tau_x -\tau_s^{-1}(\tau_s)_x]\dx;
\]
\[\mathcal{R}_3:=-\f{1}{2}\int_0^1 a'(\tau^{-2}-\tau_s^{-2})[\tau^{-1}\tau_x -\tau_s^{-1}(\tau_s)_x]\dx.
\]
Notice that $\mathcal{R}_1$ can be reformulated as
\[
\mathcal{R}_1=A \gamma \int_0^1 [\tau^{-\gamma-1}(\tau_x-(\tau_s)_x)+(\tau_s)_x(\tau^{-\gamma-1}-\tau_s^{-\gamma-1})]
\]
\[
\times [\tau^{-1}(\tau_x-(\tau_s)_x)+(\tau_s)_x(\tau^{-1}-\tau_s^{-1})]\dx.
\]
Consequently, using (\ref{m8}), (\ref{wy22}), (\ref{wy23}) and Cauchy-Schwarz's inequality, we get the estimate
\beq\label{z11}
\mathcal{R}_1 \geq C_1 \int_0^1(\tau_x-(\tau_s)_x)^2\dx-C_2 \int_0^1(\tau-\tau_s)^2\dx.
\eeq
In a similar manner, we have
\beq\label{z12}
\mathcal{R}_2 \geq -\delta_1\int_0^1(\tau_x-(\tau_s)_x)^2\dx-C_{\delta_1}\int_0^1(\tau-\tau_s)^2\dx;
\eeq
\beq\label{z13}
\mathcal{R}_3 \geq -\delta_2\int_0^1(\tau_x-(\tau_s)_x)^2\dx-C_{\delta_2}\int_0^1(\tau-\tau_s)^2\dx.
\eeq
In view of (\ref{z11})-(\ref{z13}), we choose $\delta_1,\delta_2$ to be sufficiently small and arrive at
\beq\label{z14}
-\int_0^1 (P(x,\tau)-P(x,\tau_s))_x \left(\log \left(\f{\tau}{\tau_s}\right)\right)_x\dx \geq C_3 \int_0^1(\tau_x-(\tau_s)_x)^2\dx-C_4 \int_0^1(\tau-\tau_s)^2\dx.
\eeq
It follows from (\ref{z10}), by (\ref{z14}), that
\beq\label{z15}
\f{d\mathcal{H}}{dt}+C_3 \int_0^1(\tau_x-(\tau_s)_x)^2\dx \leq C_4 \int_0^1(\tau-\tau_s)^2\dx+ \int_0^1 \varrho u_x^2 \dx,
\eeq
where we have set
\[
\mathcal{H}:=\f{\mu}{2} \int_0^1 \left(\log \left(\f{\tau}{\tau_s}\right)\right)_x^2 \dx -\int_0^1 u\left(\log \left(\f{\tau}{\tau_s}\right)\right)_x \dx.
\]
In addition, by invoking (\ref{m8}), (\ref{wy22}), (\ref{wy23}) and Cauchy-Schwarz's inequality, we see
\[
C_5 \int_0^1(\tau_x-(\tau_s)_x)^2\dx -C_6 \left(\int_0^1(\tau-\tau_s)^2\dx+\int_0^1u^2\dx\right)
\]
\beq\label{z16}
\leq \mathcal{H}\leq C_7 \left(\int_0^1(\tau_x-(\tau_s)_x)^2\dx+\int_0^1(\tau-\tau_s)^2\dx+\int_0^1u^2\dx\right).
\eeq
We multiply (\ref{z15}) both sides by $\delta_3$ and add the resulting inequality to (\ref{y10}), after choosing $\delta_3$ to be sufficiently small, to conclude that
\beq\label{z17}
\f{d}{dt}(\mathcal{E}+\delta_3 \mathcal{H})+C_8 (\mathcal{E}+\delta_3 \mathcal{H}) \leq 0.
\eeq
Then one checks easily, with the help of (\ref{y12}) and (\ref{z16}), that
\[
C_9^{-1}\left(\int_0^1(\tau_x-(\tau_s)_x)^2\dx+\int_0^1(\tau-\tau_s)^2\dx+\int_0^1u^2\dx\right)
\]
\beq\label{z18}
\leq \mathcal{E}+\delta_3 \mathcal{H} \leq C_9 \left(\int_0^1(\tau_x-(\tau_s)_x)^2\dx+\int_0^1(\tau-\tau_s)^2\dx+\int_0^1u^2\dx\right).
\eeq
Moreover, the right-hand side of (\ref{z7}) is estimated as follows:
\[
\left(\int_0^1u^2\dx+\int_0^1( P(x,\tau))_x^2\dx\right)= \left(\int_0^1u^2\dx+\int_0^1( P(x,\tau)-P(x,\tau_s))_x^2\dx\right)
\]
\[
\leq C_{10}\left(\int_0^1(\tau_x-(\tau_s)_x)^2\dx+\int_0^1(\tau-\tau_s)^2\dx\right),
\]
due to (\ref{m8}), (\ref{wy22}) and (\ref{wy23}). Therefore, (\ref{z7}) implies
\beq\label{z19}
\f{d}{dt} \int_0^1u_x^2 \dx+C_3\int_0^1 u_{xx}^2\dx \leq C_{11}\left(\int_0^1(\tau_x-(\tau_s)_x)^2\dx+\int_0^1(\tau-\tau_s)^2\dx\right).
\eeq
Again as we have done previously, multiplying (\ref{z19}) both sides by $\delta_4$ and adding the resulting inequality to (\ref{z17}), after choosing $\delta_4$ to be sufficiently small, gives
\beq\label{z20}
\f{d}{dt}\left(\mathcal{E}+\delta_3 \mathcal{H}+\delta_4 \int_0^1u_x^2\dx \right)+C_{12} \left(\mathcal{E}+\delta_3 \mathcal{H}+\delta_4 \int_0^1u_x^2\dx \right) \leq 0.
\eeq
This establishes the decay estimate
\beq\label{z21}
\|(\tau-\tau_s)(t)\|_{H^1}+\|u(t)\|_{H^1}\leq C_{13}\exp (-C_{14} t),\text{ for any }t\geq 0,
\eeq
by integrating (\ref{z20}) and noting the simple fact that
\[
C_{15}^{-1} \left(\int_0^1(\tau_x-(\tau_s)_x)^2\dx+\int_0^1(\tau-\tau_s)^2\dx+\int_0^1u^2\dx+\int_0^1u_x^2\dx\right)
\]
\[
\leq \left(\mathcal{E}+\delta_3 \mathcal{H}+\delta_4 \int_0^1u_x^2\dx \right)
\]
\[
\leq C_{15} \left(\int_0^1(\tau_x-(\tau_s)_x)^2\dx+\int_0^1(\tau-\tau_s)^2\dx+\int_0^1u^2\dx+\int_0^1u_x^2\dx\right).
\]
Finally, the decay estimate
\beq\label{z22}
\|(b-b_s)(t)\|_{H^1}\leq C_{16}\exp (-C_{17} t),\text{ for any }t\geq 0,
\eeq
is a direct consequence of (\ref{m8}), (\ref{wy22}), (\ref{z4}), (\ref{z21}) and Sobolev's inequality. The proof of Theorem \ref{sl2} is thus finished by adding (\ref{z22}) to (\ref{z21}).

\vskip0.5cm

\centerline{\bf Acknowledgement}
The research of Yang Li and Yongzhong Sun is supported by NSF of China under Grant No. 11571167 and PAPD of Jiangsu Higher Education Institutions.


\end{document}